\documentstyle[11pt,psfig]{article}

\input{psfig.sty}
 
\newtheorem{theorem}{Theorem}[section]

\newtheorem{lemma}[theorem]{Lemma}

\newtheorem{remark}[theorem]{Remark}

\newtheorem{definition}[theorem]{Definition}

\newenvironment{proof}{{\it Proof:\/}}{$\Box$\vskip 0.08in}

\begin{document}
\begin{center}
\begin{LARGE}
\baselineskip=10pt 
{Unexpected connections between Burnside Groups and Knot Theory}
\end{LARGE}
%\centerline{Version: November, 2002 - June 30, 2003}
\end{center}
% send to PNAS
%put on WEB Sept, 7, 2001
\begin{center} Mieczys{\l}aw K.~D{\c a}bkowski
and J\'ozef H.~Przytycki
\end{center}
\vspace{27pt}

\centerline {Abstract}

{\footnotesize{In Classical Knot Theory and in the  new Theory of 
Quantum Invariants substantial effort was directed toward
the search for unknotting moves on links.
We solve, in this note, several classical problems concerning unknotting moves.
Our approach uses a new concept, Burnside groups of links, 
which establishes unexpected relationship  between Knot Theory 
and Group Theory.
Our method has the potential to be used in computational biology 
in the  analysis of  DNA via tangle embedding theory, as developed by 
D.W.Sumners.  
}} 
\ \\ 
\ \\ 
Connections between knot theory and group theory
can be traced back to
Listing's pioneering paper of 1847 \cite{Lis}, in which
he considered knots, and groups of signed permutations. 
The first, well-established instance of such a connection
was provided by W.Wirtinger and M.Dehn \cite{De}. They
applied the Poincar\'e's  fundamental group of a knot
exterior to study knots and their symmetries.
The connection we describe in this note is,
on the one hand, deeply rooted in Poincar\'e's tradition and 
on the other hand it is novel and unexpected. It was discovered in
our study of the cubic skein modules of the 3-sphere and it led us to
the solution of the twenty year old Montesinos-Nakanishi conjecture.

We outline our main ideas and proofs. The complete 
exposition of this theory and its applications will be the 
subject of the sequel paper, \cite{D-P-2}.

\ \
 
\ \ 
 
\section{Open problems}\label{1} 
 
Every link can be simplified to a trivial link 
by crossing changes ( 
\psfig{figure=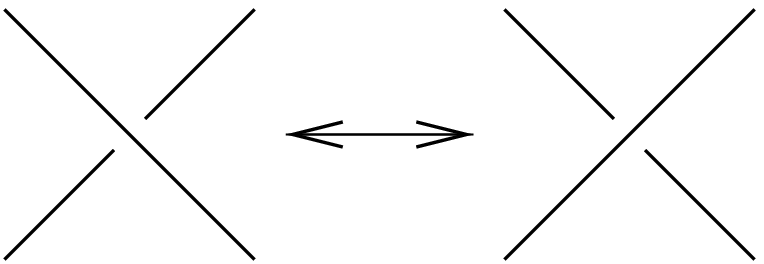,height=0.7cm}
). This observation
led to many significant developments, in particular, the construction
of the Jones polynomial of links
and the Reshetikhin-Turaev invariants of 3-manifolds.
These invariants had a great impact on modern
Knot Theory. The natural generalization of a crossing change,
 which is addressed in this paper, 
is a tangle replacement move, that is, a local modification of
a link, $L$, in which a tangle $T_1$ is replaced by a tangle $T_2$.
Several questions have been asked about which families of tangle moves
are unknotting operations.

One such family of moves, which is significant not only in Knot Theory
but also in Computational Biology, is the family of rational moves \cite{Sum}.
In this note we devote our attention to special classical cases.
 
\begin{enumerate}
\item[(1)] Nakanishi 4-move conjecture, 1979. \\
Every knot is $4$-move equivalent to the trivial knot.
\item[(2)] Montesinos-Nakanishi 3-move conjecture, 1981.\\
Every link is $3$-move equivalent to a trivial link.
\item[(3)]  Kawauchi's question, 1985.\\
Are link-homotopic links 4-move equivalent?
\item[(4)] Harikae-Nakanishi conjecture, 1992.\\
Every link is $(2,2)$-move equivalent to a trivial link.
\item[(5)] $(2,3)$-move question, 1995.\\
Is every link $(2,3)$-move equivalent to a trivial link?
\end{enumerate}

The method of Burnside groups, which we introduce, allows us to answer
questions (2), (3), (4) and (5) although Conjecture (1)
remains open.\ \\ 
We generalize questions (2), (4) and (5) to the following question:\ \\ 
\begin{enumerate}
\item[(6)] Rational moves question.\\
Is it possible to reduce every link to a trivial link
by rational $\frac{p}{q}$-moves where $p$ is a fixed prime and
$q$ is an arbitrary nonzero integer? See Def. 1.1.
\end{enumerate}
To approach our problems we define new invariants of links and call them
 the Burnside groups of links. These invariants are shown to be 
unchanged 
by $\frac{p}{q}$-rational moves.
The strength of our method lies in the fact that we are able to use
the well-developed theory of classical Burnside groups and
their associated Lie rings \cite{Mag}.  
We first describe, in more detail, how our method is applied to 
rational moves.  
In particular, we settle the Montesinos-Nakanishi and Harikae-Nakanishi
conjectures. Later we answer Kawauchi's question in detail (Section 3). 

\begin{definition}\label{1.1} 
A {\it rational $\frac{p}{q}$-move}\footnote{This move was first
considered by J.M.Montesinos \cite{Mon}.} refers to changing a link by 
replacing an identity tangle in it by a rational $\frac{p}{q}$-tangle 
of Conway (Fig.1.1).
\end{definition} 
\ \\
\centerline{\psfig{figure=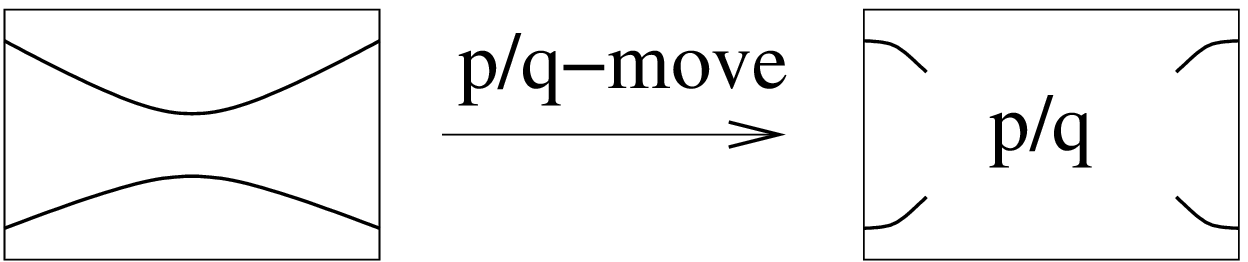,height=2.6cm}}
%\centerline{\psfig{figure=13-5-move-col.eps,height=4.6cm}}
\centerline{Fig. 1.1}

The tangles shown in Figure 1.2 are called rational tangles
and denoted by $T(a_1,a_2,...,a_n)$ in Conway's notation.
A rational tangle is the
$\frac{p}{q}$-tangle if $\frac{p}{q} =$
$a_n + \frac{1}{a_{n-1}+...+\frac{1}{a_1}}$.
Conway proved that two rational tangles are ambient isotopic
(with boundary fixed) if and only if their slopes
are equal (compare \cite{Kaw}).\\
\ \\ 
\centerline{\psfig{figure=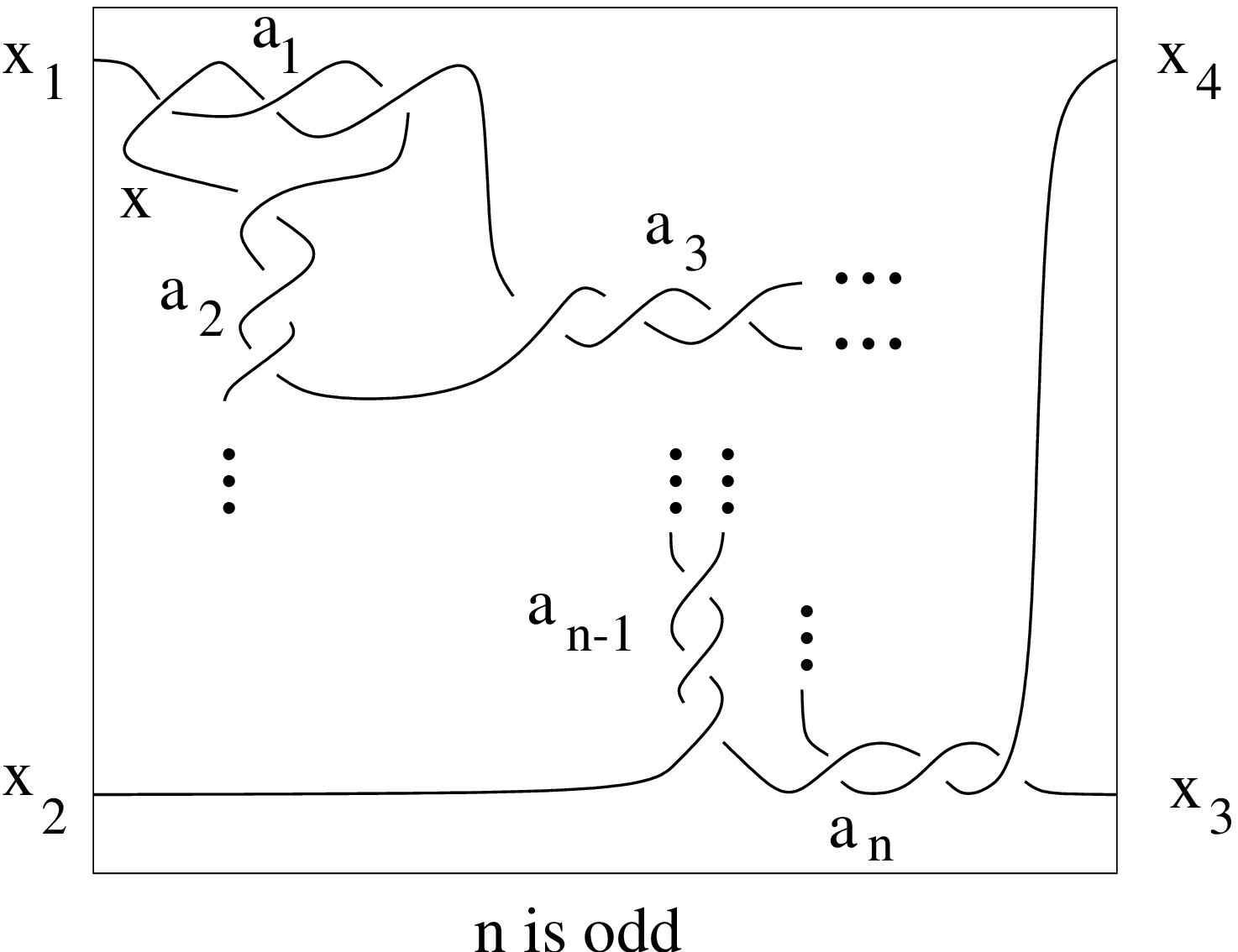,height=5.3cm}\ \ \
%\centerline{\psfig{figure=RTodd-M.eps,height=4.9cm}\ \ \  
\psfig{figure=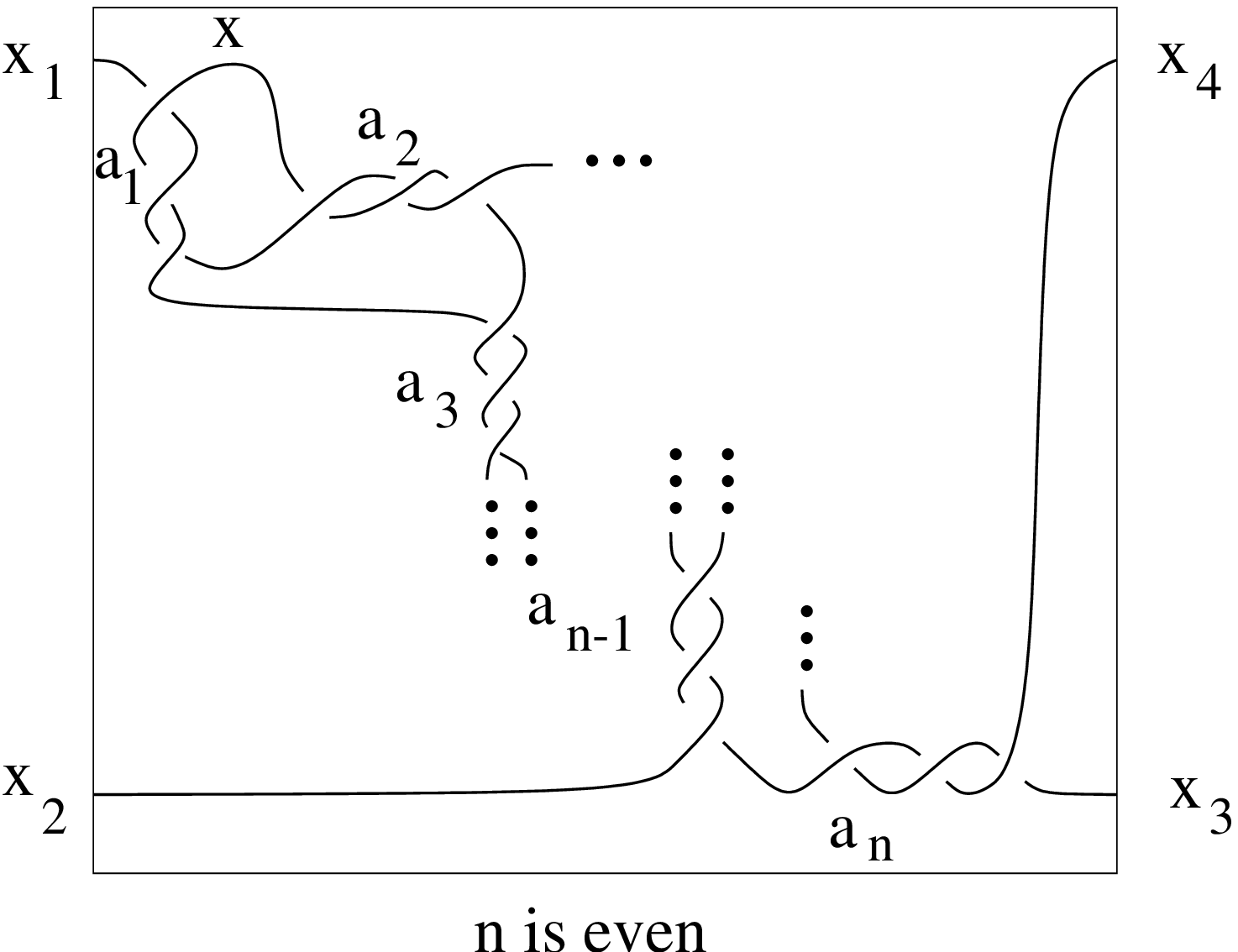,height=5.3cm}} 
\centerline{Fig. 1.2}

Rational tangles can also be viewed as tangles, which are obtained by  
applying a finite number of consecutive  
twists of the neighboring endpoints to the elementary tangle $[0]$.
\\ 
\ \\ 
Rational $\frac{p}{q}$-tangles were used  by Sumners and Ernst in  
their mathematical model of DNA recombination \cite{E-S}.  
This is a very promising development in Computational Biology.

\ \\ 
The following definition establishes a connection between two classical  
theories, Knot Theory and the Theory of Burnside 
Groups\footnote{These groups were first considered by W.Burnside in 1902,
when he asked when the group, in modern notation $B(r,n)$, is finite. 
Here  $B(r,n)$ is the quotient group of the free group on $r$ generators 
modulo the subgroup  generated by all words of the form $w^n$. It was shown 
that $B(r,n)$ is finite for $n=2,3,4,6$. On the other hand,
Novikov and Adjan proved that $B(r,n)$ is infinite for $r\geq 2$ and 
$n$ odd and sufficiently large \cite{N-A}. It is an open problem 
whether the group $B(2,5)$ is infinite, as most of experts predict,
or finite, in which case it would have $5^{34}$ elements.}. Burnside  
groups of a link play a  crucial role in our research and  
they can contribute significantly to the 
applications of Knot Theory in Computational Biology.
\ \
%Core groups, Unreduced Burnside groups.\\ 
\begin{definition}\label{1.2} 
Let $D$ be a diagram of a link $L$. We define
the associated core group $\Pi^{(2)}_D$ of $D$ by the following presentation:
generators of $\Pi^{(2)}_D$ correspond to arcs of the diagram. 
Any crossing $v_s$ yields 
the relation $r_s=y_iy_j^{-1}y_iy_k^{-1}$ where $y_i$ corresponds
to the overcrossing and $y_j,y_k$ correspond to the undercrossings 
at $v_s$ (see Fig.1.3).
\end{definition}

\begin{remark}\label{1.3}
In the above presentation of $\Pi^{(2)}_L$ one relation may
be dropped since it is a consequence of others.
\end{remark}
\ \\
\ \\
%\vspace*{1.4in} 
\centerline{\psfig{figure=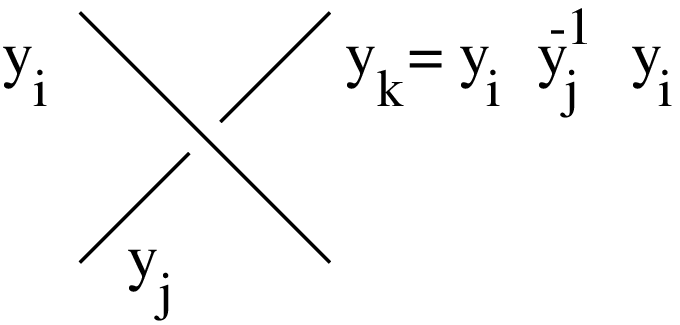,height=2.6cm}}
\centerline{Fig. 1.3}

\begin{definition}\label{1.4}
\begin{enumerate}
\item [(1)] The $n$th Burnside group of a link is the quotient of the
fundamental group of the double branched cover of $S^3$ with the link as
the branch set by its subgroup that is generated by all relations of 
the form $w^n=1$. Succinctly: $B_L(n)=\pi_1(M_L^{(2)})/(w^n)$.
\item [(2)]
The unreduced $n$th Burnside group of the unoriented link $L$ is
the quotient group 
${\hat B}_L(n) = \Pi_L^{(2)}/(w^n)$, where $ \Pi_L^{(2)}$
is the associated core group of $L$.
\end{enumerate} 
\end{definition}
 
The relation to the fundamental group of a double branched cover, 
mentioned before, is formulated below. An elementary proof
using only Wirtinger presentation and also valid for tangles 
is presented in \cite{Pr}.
 
\begin{theorem}[Wada]\label{1.5}\ 
If $D$ is a diagram of a link (or a tangle) $L$, then  
$\Pi^{(2)}_D = \pi_1(M_L^{(2)})\ast Z.$\ 
Furthermore, if we put $y_i=1$ for any fixed generator, then
$\Pi^{(2)}_D$ reduces to $\pi_1(M_L^{(2)})$. 
\end{theorem} 

The next theorem allows us to use Burnside groups to analyze 
elementary moves on links. 
 
\begin{theorem}\label{1.6} \ The groups 
$B_L(n)$ and $\hat B_L(n)$ are preserved by rational $\frac{n}{q}$-moves. 
In particular the $n$th Burnside group is preserved by $n$-moves.
\end{theorem} 
 
\begin{proof} Let $L'$ be obtained from $L$ by a $\frac{n}{q}$-move.
Then $M_{L'}^{(2)}$ is obtained from $M_L^{(2)}$
by performing the $\frac{q}{n}$-surgery. 
Such a surgery can be easily proved to preserve the
$n$th Burnside group of the fundamental group of the manifold.
\end{proof} 
 
\section{Reductions by rational moves}\label{2} 
 
We show that the answer to the rational move question, Problem (6) is negative.  
 
\begin{theorem}\label{2.1} 
\begin{enumerate} 
\item[(1)] 
The closure, $\hat \Delta_{3}^{4}$, 
 of the 3-braid $\Delta_{3}^{4}=(\sigma_1\sigma_2)^6$ 
(Figure 2.1) is not
$\frac{p}{q}$-move reducible to a trivial link for any prime number $p \geq 5$. 
\item[(2)] 
The closure of the 5-braid $\Delta_{5}^{4}=(\sigma_1\sigma_2\sigma_3\sigma_4)^{10}$
(Figure 2.2) is not 3-move reducible to a trivial link.
\end{enumerate} 
\end{theorem}  
 
%\begin{proof}
{\it Sketch of the proof:} \
We use Sanov's theorem about the structure of the Lie algebra  
associated to the Burnside group of prime exponent $p\geq 5$ 
\cite{San,V-Lee}. For $p=3$ we observe that the third Burnside group 
$B_{\hat \Delta_{5}^{4}}(3)$ is the quotient of the free Burnside group,
$B(4,3)$, by  the normal subgroup generated by relations $Q_ix_i^{-1}$, 
where the words $Q_i$ can be computed from the Fig.2.2 using the 
core relations (Fig.1.3)\footnote{One gets 
$Q_i=x_1x_2^{-1}x_3x_4^{-1}x_5x_1^{-1}x_2x_3^{-1}x_4x_5^{-1}x_i
x_5^{-1}x_4x_3^{-1}x_2x_1^{-1}x_5x_4^{-1}x_3x_2^{-1}x_1$ (see \cite{D-P-1}).}. 
To see that relations are non-trivial in  the free Burnside group 
we use the theorem of
Levi and van der Waerden about the structure of 
the associated Lie rings of Burnside groups of exponent 3
\cite{L-W}.  
%\begin{proof} \end{proof} 
 
\begin{figure}[h] 
\centerline{\psfig{figure=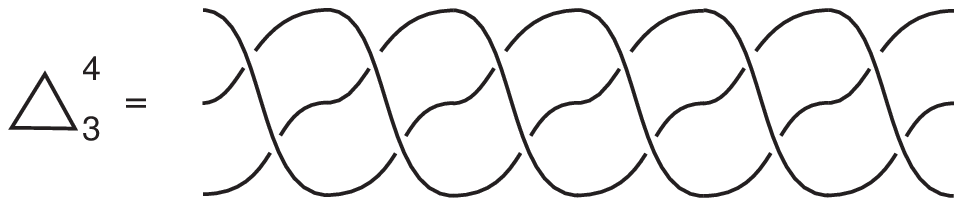,height=1.5cm}}
\centerline{Figure 2.1} 
\end{figure}
 
\begin{figure}[h] 
\centerline{\psfig{figure=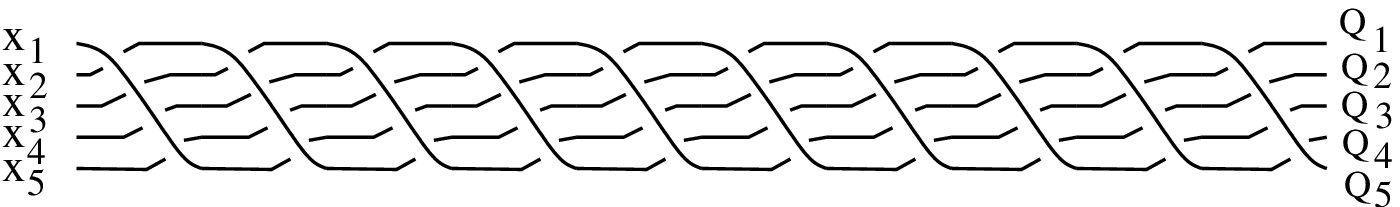,height=2.5cm}}
\centerline{Figure 2.2} 
\end{figure} 

The closed braid $\hat \Delta_{5}^{4}$ is 3-move equivalent to a link 
of 20 crossings (closure of the 5-string braid 
$(\sigma_3\sigma_4^{-1}\sigma_3\sigma_1^{-1}\sigma_2)^4$). 
It is still an open problem whether the Montesinos-Nakanishi 
conjecture holds for links up to 19 crossings. Q.Chen proved that it holds 
for links up to 12 crossings \cite{Ch}. 

The negative answers to Problems (2), (3) and (5) follow from Theorem 2.1 
for $p=3,5$ and $7$, respectively (see \cite{D-P-2} for details).

\section{Kawauchi's Question on 4-moves} 
 
In this section, we use the $4$th Burnside group of links to show that 
there is an obstruction to 4-move
reducibility of links which are link homotopically trivial.
Therefore the answer to Kawauchi's question is negative.
  
Let ${\cal W}$ denote the ``half" 2-cabling of the Whitehead link, 
which is link homotopy equivalent to the
trivial link of 3 components (Figure 3.1).\\ 
\ \\ 
\centerline{\psfig{figure=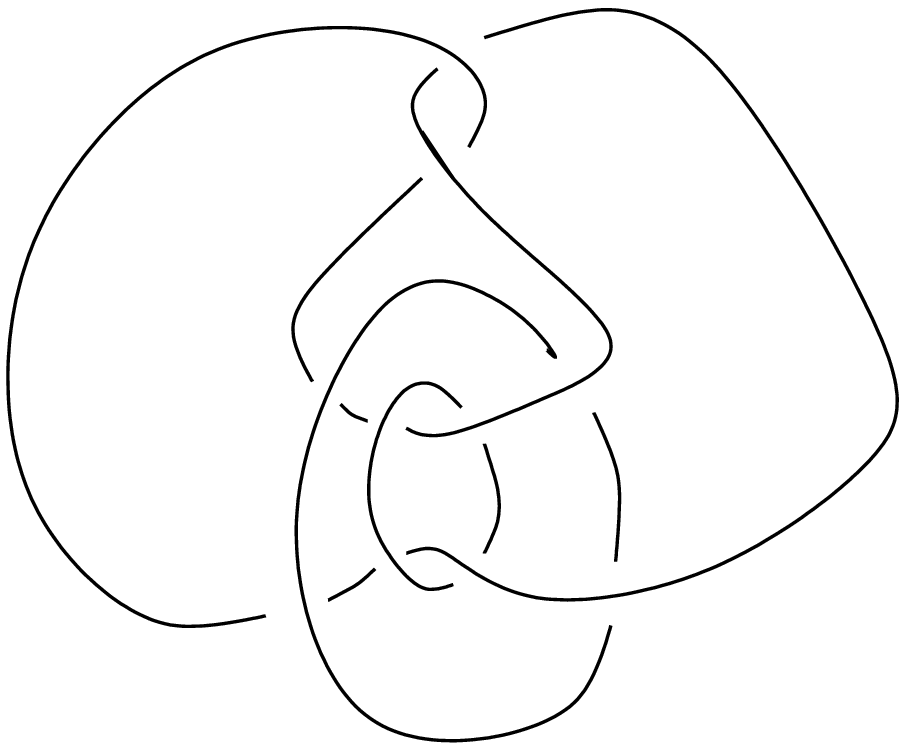,height=3.5cm}} 
\centerline{Figure 3.1} 
 
\begin{theorem}\label{3.1} 
The link ${\cal W}$ is not 4-move equivalent to a trivial link. 
\end{theorem} 
 
Our proof uses the obstruction in the Burnside group $B(2,4)$. 
The obstruction lies in the last nontrivial term of the lower central series, 
$\gamma_5$, of  
$B(2,4)$ or equivalently in the fifth term of the associated 
graded Lie ring of $B(2,4)$.  
In fact we have:  
 
\begin{lemma}\label{3.2} 
$B_{{\cal W}}(4) = B(2,4)/\gamma_5$. 
\end{lemma} 
\ \\ 
\centerline{\psfig{figure=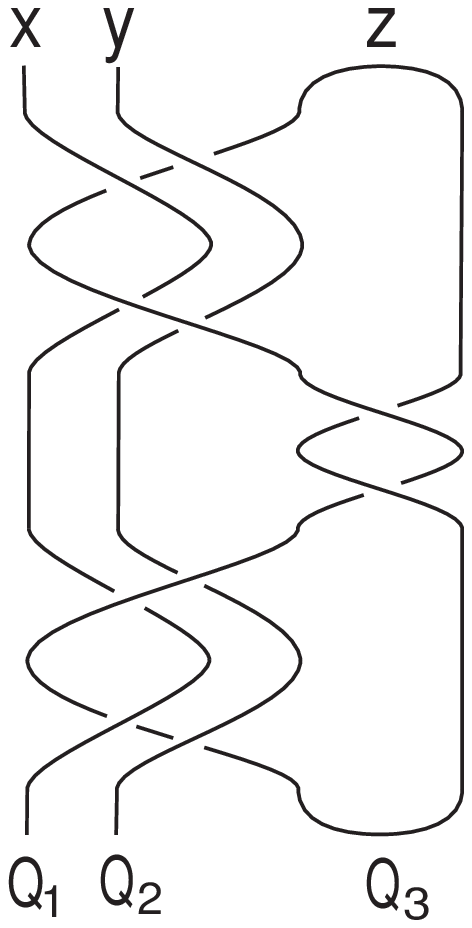,height=4.8cm}}
%\vspace*{1.6in}
\begin{center} 
Figure 3.2 Computation of relations of $B_{\cal W}(4)$.
\end{center} 
\ \\
\ \\ 
%\begin{proof}
{\ Sketch of the proof:} \ 
The lower central series of $B(2,4)$ is known to be of class 5, with the  
last term $\gamma_5=L_5$ isomorphic to $Z_2 \oplus Z_2$, \cite{V-Lee}. 
First we compute  $\pi_{1}(M_{{\cal W}}^{(2)})$ 
using a presentation of the associated core group of ${\cal W}$ by 
putting the generator $z=1$. From the diagram we obtain the relations: 
$Q_{1}x^{-1}$ and $Q_{2}y^{-1}$ (Fig.3.2), and further 
the equivalent relations \\
$R_{1}=xy^{-2}x^{2}y^{-2}x^{3}y^{-2}x^{2}y^{-2}$ and \
$R_{2}=yxy^{-2}x^{2}y^{-2}xyx^{-2}y^{2}x^{-2}$.\\ 
Therefore
$\pi_{1}(M_{{\cal W}}^{(2)})=\{x,y |\ R_{1},R_{2} \}$. \\ 
In $B_{{\cal W}}(4)$ this relations can be reduced into: \\ 
$R_{1}=x^{-1}(x^{2}y^{-2}x^{2}y^{-2})x(x^{2}y^{-2}x^{2}y^{-2})=
x^{-1}(x^{2}y^{-2})^{4}xx^{-1}(y^{2}x^{-2})^2
x(x^{2}y^{-2})^{2}$\\ 
$=x^{-1}(x^{2}y^{-2})^{4}x[x,(x^{2}y^{-2})^{2}]=
[x,(x^{2}y^{-2})^{2}]=[x,[x^{2 },y^{2}]] =
[x,y,x,y,x] \in \gamma_{5}$
(see \cite{V-Lee}). Analogously, \ $R_{2}=[y,x,y,x,y] \in \gamma_{5}$.\\
Those elements form a basis of $\gamma_{5}= Z_2\oplus Z_2$ 
considered as a $Z_2$ linear space. We verified this
fact using programs GAP, Magnus, and Magma. 
These calculations can be done manually, however they strongly depend on 
the unpublished Ph.D. thesis of
J.J. Tobin \cite{Tob}, as it was pointed out  
to us by M.Vaughan-Lee \cite{V-Lee-2}. 
 
It follows from Lemma 3.2 and the fact that $\gamma_5$ is in the center 
of $B(2,4)$, that  
$|B_{{\cal W}}(4)|=2^{10}$. On the other hand the 
abelianization of $B_{{\cal W}}(4)$ is isomorphic to $Z_{4}  
\oplus Z_{4}$. So if ${\cal W}$ were 4-reducible 
into a trivial link then the trivial link would have $3$  
components. But $B_{T_{3}}(4)=B(2,4)$ has $2^{12}$ elements, 
as predicted by Burnside and verified by Tobin  
\cite{Bu,Tob,V-Lee}. This completes the proof.  
 
%\centerline{\psfig{figure=.eps,height=1.4cm}}  
%\vspace*{1.6in}  \begin{center}  Fig. 3.3;  \end{center} 

\begin{remark}\label{3.3} 
Nakanishi showed, using Alexander modules, that 
the Borromean rings, BR, cannot be reduced to a  
trivial link by 4-moves.  
Our Burnside obstruction method also works in this case. Knowing  
that $|B(2,4)| = 2^{12}$ we can verify that $|B_{BR}(4)| = 2^{5}$. 
In addition, we conclude that ${\cal W}$ and $BR$ 
are not 4-move equivalent.
\end{remark} 
 
By Coxeter's theorem \cite{Cox} the quotient group $B_{3}/(\sigma_{i}^{4})$ 
is finite.  Therefore, for closed 3-braids, we can list all 
possible $4$th Burnside groups. This allows us to find all 4-move 
equivalence classes of closed $3$-braids \cite{D-P-2}.

\section{Limitations of the Burnside group invariant}\label{4} 
 
The method based on the Burnside group invariant has been quite successful
in the study of the unknotting property of
several
classes of tangle replacement moves. As we saw in Sections 2 and 3,
for any fixed prime number $p\geq 3$ and an arbitrary
nonzero integer $q$, rational $\frac{p}{q}$\,-\,moves are not 
unknotting operations. 
However our method has its limitations. 
We have been unable to find, by our method, obstructions 
for $\frac{n}{q}$-reduction of a link $L$, to a trivial link 
if the abelianization of the $n$th Burnside group of the link, 
$(B_{L}(n))^{(ab)} = H_{1}(M_{L}^{(2)},\mathbf{Z}_{n})$,
is a cyclic group $\textup{(i.e. }\{1\}\, or \,\mathbf{Z}_{n}\textup{)}$.
This is explained in Theorem 4.1.

Define the 
{\it restricted Burnside group of a link}, $R_{L}(n)$, 
as the quotient group $B_{L}(n)/N$,
where $N$ is the intersection of all normal
subgroups of $B_{L}(n)$ of finite indexes.  
 
\begin{theorem} Let $n$ be a power of a prime number.
Assume that the abelianization of $B_{L}(n)$ 
is a cyclic group. 
Then the restricted Burnside group, $R_{L}(n)$, is isomorphic to 
$H_{1}(M_{L}^{(2)},\mathbf{Z}_{n})$. In particular, if $B_{L}(n)$ is finite
$\textup{(e.g.\ for}\  n = 2,\,3,\,4\textup{)}$, then $B_{L}(n)$ is a cyclic 
group.
\end{theorem}

\begin{proof}
It was proved by E. Zelmanov  that 
$R(r,n)$ is finite for any $n$.
It follows that $R_{L}(n)$ is finite for any 
$n$. \\
Let $\gamma_{1}\geq\gamma_{2}\geq\cdots\geq\gamma_{i}\geq\cdots$ be the 
lower central series of $B_{L}(n)$. If
$(B_{L}(n))^{(ab)}$ (and therefore also $(R_{L}(n))^{(ab)}$)
is a cyclic group then $\gamma_{2} = \gamma_{3} = \gamma_{4}\dots $. 
Since $R_{L}(n)$ is a finite nilpotent group, therefore 
$R_{L}(n) = (B_{L}(n))^{(ab)}$. \\ 
If $B_{L}(n)$ is finite (as it is in the case of $n = 2,\,3,\,4$), 
then $R_{L}(n) = B_{L}(n)$, so $B_{L}(n)$ is a cyclic
group. 
\end{proof}
 
Therefore, the method based on the Burnside group invariant will not 
produce any obstructions for the Nakanishi 4-move conjecture and 
the Kawauchi 4-move question for a link of two components.

\section{Conclusion}\label{5} 
 
Our interest in the analysis of rational moves on links  
was inspired by our long pursuit of a program for understanding  
a 3-dimensional manifold by the knot theory which it supports.  
The method we introduced, the Burnside group of links, not only  
settles classical conjectures (e.g., the Montesinos-Nakanishi conjecture)  
but also has clear potential to be used in computational biology  
in an analysis of DNA, its recombination, action of topo-isomers, and  
analysis of protein folding and protein evolution.

%Current address:\\
 \ \\
Department of Mathematics\\
George Washington University \\
Washington, DC 20052 \\
e-mails: mdab@gwu.edu \\
przytyck@gwu.edu
\end{document}